\documentclass[11pt]{amsart}
 
\usepackage{amsfonts}
\usepackage{amssymb}
\usepackage{amscd}
\input xy
\xyoption{all}
 
\newcommand{\marginlabel}[1]%
  {\mbox{}\marginpar{\raggedleft\hspace{0pt}\bfseries\sf#1}}

\def\PP{{\textbf P}}
\def\OO{{\mathcal O}}

\def\F{\mathcal{F}}
\def\E{\mathcal{E}}
\def\G{\mathcal{G}}
\def\H{\mathcal{H}}

\def\cM{\mathcal{M}}

\def\Pic0{{\rm Pic}^0(X)}

\theoremstyle{plain}
\newtheorem*{introtheorem}{Theorem}
\newtheorem*{introcorollary}{Corollary}
\newtheorem*{introconjecture}{Conjecture}
\newtheorem{theorem}{Theorem}[section]
\newtheorem{proposition}[theorem]{Proposition}
\newtheorem{corollary}[theorem]{Corollary}
\newtheorem{lemma}[theorem]{Lemma}

\theoremstyle{definition}
\newtheorem{definition}[theorem]{Definition}
\newtheorem{remark}[theorem]{Remark}
\newtheorem{example}[theorem]{Example}
   
\newtheorem{conjecture}[theorem]{Conjecture}
\newtheorem{conjecture/question}[theorem]{Conjecture/Question}
  
\newtheorem{remark/definition}[theorem]{Remark/Definition} 
\newtheorem{terminology/notation}[theorem]{Terminology/Notation}
 
\theoremstyle{remark}

\begin{document}
 
\title{\bf Regularity on abelian varieties II: basic results on linear series 
and defining equations}
 
\author[G. Pareschi]{Giuseppe Pareschi} 
\address{Dipartamento di Matematica, Universit\`a di Roma, Tor Vergata, V.le della 
Ricerca Scientifica, I-00133 Roma, Italy}
\email{{\tt pareschi@mat.uniroma2.it}}

\author[M. Popa]{Mihnea Popa$^1$}
\footnotetext[1]{The second author was partially supported 
by a Clay Mathematics Institute Liftoff Fellowship 
during the preparation of this paper.}
\address{Department of Mathematics, Harvard University,
One Oxford Street, Cambridge, MA 02138, USA} 
\email{{\tt mpopa@math.harvard.edu}}

\maketitle


\markboth{G. PARESCHI and M. POPA}
{REGULARITY ON ABELIAN VARIETIES II}

\section*{\bf Abstract}

We apply the theory of M-regularity developed 
in \cite{us} to the study of linear series given by multiples of
ample line bundles on abelian varieties. We define an invariant of
a line bundle, called M-regularity index, which is seen to govern the
higher order properties and (partly conjecturally) the defining
equations of such embeddings. We prove a general result on the
behavior of the defining equations and higher syzygies in embeddings given by multiples of ample
bundles whose base locus has no fixed components, extending a
conjecture of Lazarsfeld proved in \cite{pareschi}. This approach also unifies
essentially all the previously known
results in this area, and is based on 
Fourier-Mukai techniques rather than representations of theta groups. 

\section{\bf Introduction}

This paper is mainly concerned with applying the theory of 
Mukai regularity (or $M$-regularity) introduced in \cite{us} to
the study of linear series given by multiples of ample line 
bundles on abelian varieties. We show that this regularity notion allows one 
to define a new invariant of a line bundle, called $M$-regularity index, which will 
be seen to roughly measure how much better one can do, given a fixed line bundle, compared 
to the standard results of the theory. Based on the main result of \cite{us} ($M$-regularity
criterion) and a related result proved here (W.I.T. regularity criterion), we show that 
all known results on such linear series can be recovered, and indeed generalized, under the 
same heading of $M$-regularity.

To make this precise, we start by recalling most of the basic results on 
ample line bundles existing
in the literature. For simplicity we state them for powers of one line bundle, 
although most hold for suitable products of possibly distinct ones. 

\begin{introtheorem}
Let $A$ be an ample line bundle on an abelian variety $X$. The following hold:
\newline
\noindent
(1) $A^2$ is globally generated.
\newline
\noindent
(2)(Lefschetz Theorem) $A^3$ is very ample.
\newline
\noindent
(3)(Ohbuchi's Theorem \cite{ohbuchi1}) If $A$ has no base divisor, then $A^2$ is very ample.
\newline
\noindent
(4)(Bauer-Szemberg Theorem \cite{bauer}) $A^{k+2}$ is $k$-jet ample, and the same holds for 
$A^{k+1}$ if $A$ has no base divisor (extending (1), (2) and (3)).
\newline
\noindent
(5)(Koizumi's Theorem \cite{koizumi}) $A^3$ gives a projectively normal embedding.
\newline
\noindent
(6)(Ohbuchi's Theorem \cite{ohbuchi2}) $A^2$ gives a projectively normal embedding if and only 
if $0_X$ does not belong to a finite union of translates of the base locus of $A$
(cf. \S5 for the concrete statement).
\newline
\noindent
(7)(Mumford's Theorem \cite{mumford1}, \cite{kempf1}) For $k\geq 4$, 
the ideal of $X$ in the embedding given 
be $A^k$ is generated by quadrics. In the embedding given by $A^3$ it is generated 
by quadrics and cubics.
\newline
\noindent
(8)(Lazarsfeld's Conjecture \cite{pareschi}, extending results of Kempf \cite{kempf2})
$A^{p+3}$ satisfies property $N_p$
(extending (5) and (7)).
\newline
\noindent
(9)(Khaled's Theorem \cite{khaled}) If $A$ is globally generated, then the ideal of $X$ in the 
embedding given by $A^2$ is generated by quadrics and cubics.
\end{introtheorem}

These results turn out to be -- some quick while others non-trivial -- consequences 
of the general global generation criterion in \cite{us}, called the $M$-regularity
criterion. Together with a more technical extension (the W.I.T. regularity 
criterion), described below, this approach yields new results and extensions 
as well. To introduce them, we first need some terminology.

Let $X$ be an abelian variety of dimension $g$ over an algebraically closed 
field, with dual abelian variety $\hat{X}$, and let
${\mathcal P}$ be a suitably normalized Poincar\'e line bundle on $X\times \hat X$. 
The Fourier-Mukai functor \cite{mukai} is the derived functor associated to the 
functor $\hat{\mathcal
S}(\F)={p_{\hat X}}_*(p_X^*\F\otimes {\mathcal P})$ from ${\rm Mod}(X)$ to 
${\rm Mod}({\hat X})$. A sheaf $\F$ on
$X$ is said to satisfy the Weak Index Theorem (W.I.T.) with index $i(\F)=k$ if
$R^i\hat{\mathcal S}(\F)=0$ for all $i\ne k$, in which case
$R^k\hat{\mathcal S}(\F)$ is simply denoted $\hat \F$. A weaker condition, introduced 
in \cite{us}, is the following: $\F$ is called \emph{$M$-regular} if ${\rm codim}(
{\rm Supp}~R^i\hat{\mathcal S}(\F))> i$ for all $i>0$.
Moreover, we will consider the \emph{Fourier jump locus} of $\F$ to be the locus of
$\xi\in \hat{X}$
where $h^0(F\otimes P_\xi)$ is different from the generic value (where $P_\xi$
is the line bundle on $X$ classified by $\xi$). 

\noindent
Given an ample line bundle $A$ on $X$, we define the \emph{$M$-regularity index} of $A$ 
to be
$$m(A):={\rm max}\{l~|~A\otimes 
m_{x_1}^{k_1}\otimes \ldots \otimes m_{x_p}^{k_p}~{\rm is ~}M{\rm -regular~ for 
~all~distinct~}$$
$$x_1,\ldots,x_p\in X {\rm~with~} \Sigma k_i=l\}.$$
A first result is that this invariant 
governs the higher order properties of embeddings obtained from $A$.

\begin{introtheorem}
If $A$ is an ample line bundle on $X$ and $k\geq m(A)$, 
then $A^{\otimes(k+2-m(A))}$ is $k$-jet ample.
\end{introtheorem}

\noindent
It is not hard to see that for example $m(A)\geq 1$ if and only if $A$ has no base 
divisor. The theorem thus recovers and extends the results of Lefschetz, Ohbuchi 
and Bauer-Szemberg mentioned above. Most interestingly though, this shows 
that results with seemingly unrelated proofs are simply steps in a hierarchy 
of regularity conditions. It is interesting to note also that the $M$-regularity 
indices are quite intimately related to the Seshadri constants measuring 
local positivity (cf. \cite{lazarsfeld2} for the case of abelian varieties); 
we will approach this in detail somewhere else.

By reversing the natural order in the body of the paper, 
the results presented in what follows suggest that it is quite natural to expect 
that a similar phenomenon governs the 
behavior of defining equations of $X$, and more generally higher syzygies, in embeddings
of this kind.

\begin{introconjecture} Let $p\ge m$ be non-negative integers.
If $A$ is ample and $m(A)\ge m$, then $A^k$ satisfies
$N_p$ for any $k\ge p+3-m$. 
\end{introconjecture}

\noindent
This extends Lazarsfeld's conjecture, which is the statement for $m=0$, 
meaning no conditions on $A$. That case has already been proved in \cite{pareschi}, 
by methods which are included in, and provide a basis for, the strategy adopted here.

The main result of this paper is a proof, and also a strenghtening, of the conjecture 
above for $m=1$, i.e. for line bundles whose base locus has no fixed
components. We first recall the terminology introduced in \cite{green}: 
property $N_p$ for a very ample line bundle means that $I_{X,L}$, the
homogeneous ideal of $X$ in the corresponding embedding, is generated  
by quadratic forms, and also that -- up to the $p$-th step -- the
higher syzygies between these forms are generated in the lowest
possible degree, i.e. by linear ones. Thus, in this language, 
the property that $I_{X,L}$ be generated by quadrics is condition
$N_1$. Moreover, the property of "being off" by $r$ from $N_p$ was formalized 
in \cite{pareschi} into property $N_p^r$ (cf. \S6 for details).
In a word, $N_p^0$ is equivalent to $N_p$, and $N_p^1$ means 
that $I_{X,L}$ is generated by quadrics and cubics.

\begin{introtheorem}\emph{(${\rm char}(k)$ does not divide $(p+1)$ and
$(p+2)$.)} 
Let $A$ be an ample line bundle on $X$, with no base divisor. Then:
\newline
\noindent
(a) If $k\ge p+2$ then $A^{k}$ satisfies property $N_p$.
\newline
\noindent
(b) More generally, if $(r+1)(k-1)\ge p+1$ then $A^k$ satisfies
property $N_p^r$.
\end{introtheorem}

\noindent
The first instance of this theorem, worth emphasizing individually, is
the following:

\begin{introcorollary}(${\rm char}(k)\ne 2,3$.)
Let $A$ be an ample line bundle on $X$, with no base divisor. Then: 
\newline
\noindent
(a) If $k\ge 3$ then $I_{X,A^{k}}$ is generated by quadrics.
\newline
\noindent
(b) $I_{X,A^2}$ is generated by quadrics and cubics.
\end{introcorollary} 
\noindent
(Note in particular the improvement of Khaled's result above.)

For consistency reasons, we note that Ohbuchi's projective normality result (6)
does not integrate in the discussion above, and does indeed suggest what happens 
in the cases left out by the above conjecture. However, it can still be obtained in 
a similar way, and in \S5 we sketch its proof as a toy version of that of the syzygy theorem.

As previously mentioned, all the proofs of the statements above are based on the 
basic $M$-regularity theorem, which we recall below. 

\begin{introtheorem}{\bf ($M$-regularity criterion, \cite{us} Theorem 2.4.)}
Let $\F$ be a coherent sheaf and $L$ an invertible sheaf supported on a 
subvariety $Y$ of the abelian variety $X$ (possibly $X$ itself). If both $\F$ and $L$ are 
$M$-regular as sheaves on $X$, then $\F\otimes L$ is globally generated.
\end{introtheorem}

\noindent 
This has to be combined with a refined study, in a relative setting,
of the notion of skew-Pontrjagin product introduced in \cite{pareschi},
and also with a different (but related) regularity criterion, which 
we prove here following a similar strategy. The new statement 
needs stronger hypotheses on the sheaf $\F$, but provides specific information
about the loci where suitable tensor products are not globally generated.

\begin{introtheorem}{\bf (W.I.T. regularity criterion.)}
Let $A$ be an ample line bundle on $X$. Let also $F$
be a locally free sheaf on $X$ such that
\newline
\noindent
(1) the Fourier-jump locus $J(F)$ is finite. 
\newline
\noindent
(2) $F^\vee$ satisfies the W.I.T. with index $i(F^\vee)=g$.
\newline
\noindent
(3) the torsion part of $\hat {F^\vee}$ is a sum of (possibly
zero) skyscraper sheaves on the points of $J(F)$. 
\newline
\noindent
Then there is an inclusion of non-generation loci:
$$B(F\otimes A)\subset \bigcup_{\xi \in J(F)}B(A\otimes P_\xi^\vee).$$
(For a sheaf $\F$, we denote by $B(\F)$ the locus where $\F$ is not globally generated.)
\end{introtheorem}

\noindent
It is interesting to note that the W.I.T. regularity criterion applies to some sheaves
for which the $M$-regularity criterion does not apply and conversely. 

The underlying principle in this article is the use of vanishing theorems and Fourier-Mukai 
methods for vector bundles, or even arbitrary coherent sheaves, in the study of linear series. 
This completely bypasses methods based on theta-functions and representations of theta-groups 
originating in \cite{mumford1} (and employed in the original proofs of most of the 
previously known results)
as it was somewhat hinted that it could be possible by earlier work of Kempf. 
Better still, the main advantage of the present methods is that they apply to a much 
wider spectrum of problems on abelian varieties, as it is described in \cite{us}.

The paper is organized as follows: in Section 2 we recall the main terminology 
and results from \cite{us}, and we introduce further notions of generation 
of sheaves. Section 3 contains the definition of the $M$-regularity index and the 
corresponding result on higher order properties of embeddings. Sections 4 and 5 are 
devoted to a rather long list of technical results needed in the study of defining 
equations. In the former we prove the W.I.T. regularity criterion, while in the latter
we introduce the notion of relative skew Pontrjagin product and study its properties 
under various operations. Finally, Section 6 is devoted to the main results of this 
paper, on defining equations and syzygies of abelian varieties embedded with powers 
of line bundles whose base locus has no fixed components.

\noindent
\textbf{Acknowledgements.} We would like to thank R. Lazarsfeld for some 
very useful conversations.

\section{\bf Background and preliminary results}

In what follows $X$ will be an abelian variety over an algebraically closed 
ground field $k$. Restrictions on ${\rm char}(k)$ will be specified along the paper.
We denote by $\hat X$ the dual of $X$, which we identify with $\Pic0$. Given
$\xi\in \hat X$, $P_\xi$ will denote the line bundle on $X$ classified
by $\xi$. For a positive integer $n$, $X_n$ will denote the group of $n$-torsion
points of $X$. When it appears in the text, we will always be in a situation
where ${\rm char}(k)$ does not divide $n$.

\subsection*{\bf Various notions of generation of sheaves.}
Let $\F$ be an arbitrary coherent sheaf on $X$. The support
of the cokernel of the evaluation map $H^0(\F)\otimes\OO_X\rightarrow \F$
will be referred to as the \emph{non-generation locus} of $\F$ and denoted $B(\F)$. 
The usual notions of global generation and generic global generation mean
that $B(\F)$ is empty or a proper subset respectively.
On abelian varieties it is useful to consider weaker notions of generation, which 
can be in fact defined on any irregular variety:
$\F$ is said to be \emph{continuously globally generated} (cf. \cite{us} \S2)
if the map
$$\bigoplus_{\xi\in U}H^0(\F\otimes P_\xi)\otimes
P_\xi^\vee\rightarrow \F$$ 
is surjective for \emph{ any non-empty Zariski-open set $U\subset \hat X$}.
For a line bundle $A$ this just means that
$\bigcap_{\xi\in U}B(A\otimes P_\xi)$ is empty. 
In what follows we introduce an even weaker variant, needed in the sequel. 

\begin{definition}
Given a sheaf $\F$, we define its \emph{Fourier jump locus} as the 
locus $J(\F)\subset \hat X$ consisting of $\xi\in \hat X$ where
$h^0(\F\otimes P_\xi)$ jumps, i.e. it is different from its minimal value over $\Pic0$. 
$\F$ is said to be \emph{weakly continuously generated} if the map
$$\bigoplus_{\xi\in U}H^0(\F\otimes P_\xi)\otimes
P_\xi^\vee\rightarrow \F$$ 
is surjective for \emph{any  non-empty Zariski-open set $U\subset \hat X$ containing
$J(\F)$}. Continuous global generation obviously implies weak continuous generation
and the two notions are equivalent if the Fourier jump locus of $\F$ is empty. 
\end{definition}

\begin{remark}\label{generation} 
The following facts are easy to check (cf. also \cite{us} Remark 2.11):
\item{(a)} If $\F$ is weakly continuously generated, then there exist
$\xi_1,\dots,\xi_k\in \hat X$ such that the map $\oplus_{i=1}^kH^0(\F\otimes
P_{\xi_i})\otimes P_{\xi_i}^\vee\rightarrow \F$ is surjective.
\item{(b)} If $F$ is continuously globally generated then there is a positive integer $N$ such
that \emph{for  general $\xi_1,\dots ,\xi_N\in \hat X$, } the map
$\oplus_{i=1}^NH^0(\F\otimes P_\xi)\otimes P_\xi^\vee\rightarrow \F$ is
surjective.
\item{(c)} If $\F$ is weakly continuously generated and 
$J(\F)$ is finite, say $J(\F)=\{\xi_1,\dots ,\xi_n\}$, then there is a positive integer
$N$ such that \emph{for general $\xi_{n+1},\dots \xi_{n+N}\in \hat X$ }
the map $\oplus_{i=1}^{n+N}H^0(\F\otimes P_\xi)\otimes P_\xi^\vee\rightarrow \F$
is surjective.
\end{remark}

\noindent 
The following lemma, proved in \cite{us} Proposition 2.12, shows how to 
produce global generation from continuous global generation.

\begin{lemma}
If $\F$ is a continuously globally generated sheaf on $X$
and $L$ is a continuously globally generated sheaf on $X$ which is everywhere 
of rank $1$ on its support, then $\F\otimes L$ is globally generated. 
\end{lemma}
\noindent
The proposition below is a variation of this result, relating the notions of weak
continuous generation and generic global generation. At least if the
Fourier jump locus of $\F$ is {\it finite}, one can describe the non-generation locus of
$\F\otimes L$ in terms of $J(\F)$.

\begin{proposition}\label{ggg} Let $\F$ be a weakly continuously generated
sheaf on $X$.
\newline
\noindent
(a) If $\E$ is a sheaf
such that $\E\otimes P_\xi$ is generically globally generated for any
$\xi\in \hat X$, then $\F\otimes \E$ is generically globally generated.
\newline
\noindent
(b) If $\E\otimes P_\xi$ is globally generated for any $\xi\in
\hat X$, then $\F\otimes \E$ is globally generated.
\newline
\noindent
(c) Assume that the Fourier jump locus
$J(\F)$ is finite, and let $L$ be a continuously globally generated line bundle
on $X$. Then $B(\F\otimes L)\subset \bigcup_{\xi\in J(\F)}B(L\otimes P_\xi^\vee).$
\end{proposition}
\begin{proof} (a) By Remark \ref{generation}(a), the map $\oplus_{i=1}^k
H^0(\F\otimes P_{\xi_i})\otimes \E\otimes P_{\xi_i}^\vee\rightarrow \F\otimes
\E$ is surjective. Therefore we have the inclusion of non-generation loci
$B(\F\otimes \E)\subset \bigcup_{i=1}^kB(\E\otimes P^\vee_{\xi_i})$. This
proves the assertion since, by hypothesis, $B(\E\otimes
P^\vee_{\xi_i})$ are proper subvarieties. This same argument also proves (b). 
\newline
\noindent
(c) If $J(\F)=\{\xi_1,\dots ,\xi_N\}$, then  the map
$\oplus_{i=1}^{n+N}
H^0(\F\otimes P_{\xi_i})\otimes L\otimes P_{\xi_i}^\vee\rightarrow
\F\otimes L$ is surjective for general $\xi_{n+1},\dots \xi_N\in \hat X$
(cf. Remark \ref{generation}(c)). Therefore
$B(\F\otimes L)$ is contained in the union of $\cup_{i=1}^n B(L\otimes
P_{\xi_i})$ with the intersection -- for all
$\xi_{n+1},\dots ,\xi_{n+N}$ general in $\hat X$ -- of
$\cup_{i=n+1}^{n+N}B(L\otimes P_{\xi_i}^\vee)$. Since $L$ is a 
continuously globally generated \emph{line bundle}, this intersection is empty. 
This implies that $B(\F\otimes L)\subset \cup_{i=1}^nB(L\otimes P_{\xi_i}^\vee)$.
\end{proof}

\subsection*{\bf Fourier-Mukai functor, index theorems and $M$-regularity.} 

According to Mukai \cite{mukai}, one considers the left-exact functor
$\hat{\mathcal S}$ from the category of $\OO_X$-modules to the category
of $\OO_{\hat X}$-modules defined as $\hat{\mathcal S}(\F)={p_{\hat
X}}_*(p_X^*\F\otimes {\mathcal P})$. Mukai's main result \cite{mukai} 
Theorem 2.2 is that the derived functor ${\bf R}\hat {\mathcal S}$
establishes an equivalence of categories between ${\bf D}(X)$ and ${\bf D}(\hat X)$.
A sheaf $\F$ on  $X$ is said to satisfy the \emph{Index Theorem
(I.T.) with index $i(\F)=k$} if $H^j(F\otimes P_\xi)=0$ for any $\xi\in
\hat X$ and any $j\ne k$. More generally, $\F$ is said to satisfy the 
\emph{Weak Index Theorem (W.I.T.) with index $i(F)=k$} if $R^j\hat{\mathcal S}(\F)=0$
for $j\ne k$. In this case $R^{i(\F)}\hat{\mathcal S}(\F)$ is simply
denoted $\hat \F$. 
A useful consequence of Mukai's theory is the following lemma (\cite{mukai}, Cor.2.5):

\begin{lemma}\label{exts}
If $\F$ and $\G$ both satisfy W.I.T., then there is a natural isomorphism
$$\phi:{\rm Ext}^j(\F,\G)\cong {\rm Ext}^{j+i(\F)-i(\G)}(\hat \F,\hat \G).$$
\end{lemma}

\noindent
We recall from \cite{us} the following weakening of the W.I.T. condition with 
index $0$, and the main regularity result proved there.

\begin{definition}
A sheaf $\F$ on $X$ is said to be \emph{Mukai-regular}
(or simply \emph{$M$-regular}) if $R^i\hat{\mathcal S}(\F)$ 
is supported in {\rm codim}ension $> i$ for any $i>0$, where for $i=g$ this
means that the support $R^g\hat{\mathcal S}(\F)$ is empty. This happens 
in particular if the cohomological support loci 
$$V^i(\F):=\{\xi~|~h^i(\F\otimes P_{\xi})\neq 0\}\subset \Pic0$$
have codimension $>i$ for all $i$. 
\end{definition}

\begin{theorem}{\bf ($M$-regularity criterion, 
\cite{us} Theorem 2.4 and Proposition 2.13.)}\label{F-reg} 
Let $\F$ be an $M$-regular sheaf on $X$, possibly supported 
on a subvariety $Y$ of $X$. Then the following hold:
\newline
\noindent
(a) $F$ is continuously globally generated.
\newline
\noindent
(b) Let also $A$ be a line bundle on $Y$, continuously globally generated 
as a sheaf on $X$.
Then $F\otimes A$ is globally generated.
\end{theorem}

\section{\bf The $M$-regularity index and properties of embeddings}

In this section we show that the concept of $M$-regularity is well-adapted 
to the study of linear series, and provides a uniform point of view on the 
study of (higher order) properties of embeddings. 
This will serve as an introduction to the deeper facts on defining equations
treated in the subsequent sections via stronger regularity techniques.
We first need to recall the notion of $k$-jet ampleness (cf. e.g. \cite{sommese}
in general and \cite{bauer} in the context of abelian varieties).

\begin{definition}
A line bundle $A$ is called $k$-\emph{jet ample}, $k\geq 0$, if the restriction map
$$H^0(A)\longrightarrow H^0(A\otimes \OO_X/
m_{x_1}^{k_1}\otimes \ldots \otimes m_{x_p}^{k_p})$$ 
is surjective for any distinct points $x_1,\ldots,x_p$ on $X$ such that $\Sigma k_i =k+1$.
\end{definition}

\begin{remark}
In particular $0$-jet ample means 
globally generated, $1$-jet ample means very ample.
The notion of $k$-jet ampleness is stronger than a related notion of $k$-very 
ampleness, which takes into account  $0$-dimensional subschemes of length
equal to $k+1$.
\end{remark}

\begin{lemma}\label{kva}
For an ample line bundle $A$ on the abelian variety $X$, the following are equivalent:
\newline
\noindent
(i) $A$ is $k$-jet ample.
\newline
\noindent
(ii) $A\otimes 
m_{x_1}^{k_1}\otimes \ldots \otimes m_{x_p}^{k_p}$
satisfies I.T. with index $0$ for all $x_1,\ldots,x_{p}$ such that $\Sigma k_i = k+1$.
\newline
\noindent
(iii) $A\otimes 
m_{x_1}^{k_1}\otimes \ldots \otimes m_{x_l}^{k_l}$ is globally generated for all  
$x_1,\ldots,x_{l}$ such that $\Sigma k_i= k$.
\end{lemma}
\begin{proof}
This is based on the immediate fact that, since $h^1(A)=0$ as we are on an 
abelian variety, $k$-jet ampleness is equivalent to the vanishing
$$H^1(A\otimes
m_{x_1}^{k_1}\otimes \ldots \otimes m_{x_p}^{k_p})=0$$
for all $x_1,\ldots,x_{p}$ such that $\Sigma k_i = k+1$.   
If $A$ is $k$-jet ample, then so is any translate, thus (i) is equivalent to 
(ii) by the very definition. 
The equivalence with (iii) also follows quickly, since the required global 
generation is equivalent to the surjectivity of
$$H^0(A\otimes 
m_{x_1}^{k_1}\otimes \ldots \otimes m_{x_l}^{k_l})\longrightarrow 
H^0(A\otimes 
m_{x_1}^{k_1}\otimes \ldots \otimes m_{x_l}^{k_l}\otimes \OO_X/m_x) $$
for every $x\in X$.
\end{proof}

\noindent 
The key definition is given below. We note that it is suggested naturally by Theorem 
\ref{F-reg} and Lemma \ref{kva}, and as a result the theorem which follows is almost 
tautological.

\begin{definition}
The $M$-\emph{regularity index} of $A$ is defined as 
$$m(A):={\rm max}\{l~|~A\otimes 
m_{x_1}^{k_1}\otimes \ldots \otimes m_{x_p}^{k_p}~{\rm is ~}M{\rm -regular~ for 
~all~distinct~}$$
$$x_1,\ldots,x_p\in X {\rm~with~} \Sigma k_i=l\}.$$
\end{definition}

\noindent
The following description provides more intuition for this definition.

\begin{proposition}\label{int}
We say that a line bundle $A$ is $k$-jet ample in codimension $r$ if
the set of points $x$ for which there exist $x_2,\ldots, x_p$ and 
$k_1,\ldots ,k_p$ with $\Sigma k_i=k+1$ such that $h^1(A\otimes m_x^{k_1}
\otimes m_{x_2}^{k_2}\otimes \ldots \otimes m_{x_p}^{k_p})>0$ has 
codimension $r$ in $X$. We have that $m(A)\geq k+1$ if 
$A$ is $k$-jet ample in codimension $\geq 2$.
\end{proposition}
\begin{proof}
If we assume that $m(A)<k+1$, then there exist $x_1,\ldots, x_p$ and 
$k_1,\ldots, k_p$ with $\Sigma k_i=k+1$ such that the set 
$$\{y\in X~|~h^i(t_y^* A \otimes m_{x_1}^{k_1}
\otimes m_{x_2}^{k_2}\otimes \ldots \otimes m_{x_p}^{k_p})>0\}$$
has codimension $\leq 1$. Since this is the same as the set 
$$\{y\in X~|~h^i(A \otimes m_{x_1-y}^{k_1}
\otimes m_{x_2 -y}^{k_2}\otimes \ldots \otimes m_{x_p -y}^{k_p})>0\},$$
the assertion follows immediately.
\end{proof}

\begin{example}{\bf (Small values of $m(A)$.)}\label{base_div}
If $A$ is an ample line bundle, then $m(A)\geq 1$ if and only if $A$ 
does not have a base divisor. Also, if $A$ gives a birational map which 
is an isomorphism outside a codimension $2$ subset, then $m(A)\geq 2$. 
Both assertions follow immediately from the proposition above.
\end{example}

\begin{example}({\bf Abelian surfaces.})
(i) Let $A$ be a polarization of type $(1,2)$ on an abelian surface. It is an immediate 
consequence of the Decomposition Theorem (cf. \cite{lange} 4.3.1)
that $A$ has a base divisor if and only if $X$ is a product of 
elliptic curves $E$ and $F$ and $A=\OO_X(E+2F)$. On the other hand, it is not hard to 
see that we always have $m(A)\leq 1$ (for example $A\otimes m_x^2$ is not $M$-regular, 
where $x$ is any point on $X$). Thus $m(A)=1$ exactly when the pair $(X,A)$ 
is not of the above form, while otherwise it is $0$. Polarizations of type $(1,3)$
are globally generated, and so $m(A)\geq 1$, but again an argument similar
to Proposition \ref{int} shows that $m(A)=1$.
\newline
\noindent
(ii) Let $A$ be a polarization of type $(1,4)$. If $X$ is general
(cf. \cite{lange} Ch.10 \S5, or the original \cite{blvs}, 
for more precise conditions), then $A$ 
gives a birational morphism to $\PP^3$ which is not an embedding, and 
whose exceptional set is a curve (we thank the referee for pointing this out to us).
However from the properties of this curve and the fact that the map 
separates tangent vectors outside a codimension $2$ subset, it follows 
easily that $m(A)\geq 2$, although one cannot directly use Proposition \ref{int}.
On the other hand, the special such abelian variety is a cover of a product of 
elliptic curves, and in that case $m(A)=1$.
\newline
\noindent
(iii) On a general abelian surface (more precisely, by Reider's Theorem \cite{lange}
10.4.1, one on which there are no elliptic curves $C$ such that $C\cdot A=2$), 
a polarization $A$ of type $(1,d)$ with $d\geq 5$ is 
very ample, and so $m(A)\geq 2$.
\end{example}

Based on this definition we obtain the following theorem, which extends and 
places in a natural setting the basic results on embeddings given by multiples 
of ample line bundles existing in the literature.

\begin{theorem}
If $A$ and $M_1,\ldots ,M_{k+1-m(A)}$ are ample line bundles on $X$, $k\geq m(A)$, 
then $A\otimes M_1\otimes \ldots\otimes M_{k+1-m(A)}$ is $k$-jet ample. In particular
$A^{\otimes(k+2-m(A))}$ is $k$-jet ample.
\end{theorem}
\begin{proof}
By definition, $A\otimes  
m_{x_1}^{k_1}\otimes \ldots \otimes m_{x_p}^{k_p}$ is $M$-regular 
for any $x_1,\ldots,x_p\in X$ as long as $\Sigma k_i = m(A)$.
This in turn implies that 
$M_1\otimes A\otimes  
m_{x_1}^{k_1}\otimes \ldots \otimes m_{x_p}^{k_p}$ is globally generated, by Theorem \ref{F-reg}.
Now, by Lemma \ref{kva}, this is the same as saying that $M_1\otimes A\otimes 
m_{x_1}^{k_1}\otimes \ldots \otimes m_{x_l}^{k_l}$
satisfies I.T. with index $0$ for all $x_1, \ldots, x_l$ and $\Sigma k_i=m(A)+1$. 
As this is a strong form of $M$-regularity, it allows 
us to continue the same procedure inductively to obtain the desired conclusion.
\end{proof}

In particular, for small values of $m(A)$ (namely $0$ and $1$), this recovers as particular cases
the theorems of Lefschetz, Ohbuchi, and more generally Bauer-Szemberg, mentioned 
in the introduction (cf. also Example \ref{base_div}).

We conclude by noting that a deeper reason for considering these invariants 
is that they seem to link in a natural 
way the geometry of the abelian variety in the embedding given by a line bundle 
with the equations, and more generally the syzygies, of that embedding. This will be 
explained in detail at the end of \S6.

\section{\bf Cohomological criteria for weak continuous generation}

In this section we provide a criterion, based on the weak 
index theorem, for the weak continuous generation of locally free sheaves on abelian
varieties. To put the result into perspective, we recall that 
a key step in the proof of Lazarsfeld's conjecture \cite{pareschi} was based on 
the fact that if $F$ is a locally free sheaf on $X$ satisfying I.T. with index $0$, 
and $A$ is an ample line bundle, then $F\otimes A$ is globally generated. 
Theorem \ref{F-reg} above, proved in
\cite{us}, provides a generalization of that criterion widely extending its
range of applicability.
For the purposes of this paper, we also need the following different generalization of 
the result mentioned above, based on even weaker hypotheses (but note the locally
freeness assumption):

\begin{theorem}{\bf (W.I.T. regularity criterion)}\label{WIT}
Let $F$ be a locally free sheaf on $X$ such that
$F^\vee$ satisfies W.I.T. with index $g$ and the torsion part of
$\widehat{F^\vee}$ is a torsion-free sheaf on a reduced subscheme of
$X$. Then the following hold:
\newline
\noindent
(a) $F$ is weakly continously generated.
\newline
\noindent
(b) Let morever $A$ be a continuously globally generated line bundle on a
subvariety of $X$. Then 
\begin{itemize}
\item[(i)] $F\otimes A$ is generically globally generated.
\item[(ii)] If the Fourier-jump locus $J(F)$ is finite then $B(F \otimes
A)\subset \bigcup_{\xi\in J(F)}B(A\otimes P_\xi)$.
\end{itemize}
\end{theorem}

\begin{corollary}\label{bpf}
Let $F$ and $A$ be a locally free sheaf, respectively an invertible sheaf on $X$. 
If $F$ satisfies the hypotheses of Theorem \ref{WIT} and $A$ is globally 
generated, then $F\otimes A$ is globally generated.
\end{corollary}

\noindent 
The corollary follows immediately from Theorem \ref{WIT}(a) and
Proposition \ref{ggg}(b). Turning to the proof of
Theorem \ref{WIT}, note that, in view of Proposition \ref{ggg}(c), 
the only thing to prove is part(a). 
This in turn follows in a standard way (cf. \cite{pareschi} \S2(b)
or \cite{us} \S2) from the following corresponding  
generalization of a result of Mumford-Kempf-Lazarsfeld type on multiplication
maps (cf. \cite{us} Theorem 2.5).

\begin{lemma}\label{mult} 
Let $F$ be a locally free sheaf on X
satisfying the hypotheses of Theorem \ref{WIT} and let $\H$ be a coherent sheaf on
$X$ satisfying I.T. with index $0$.
Then the sum of multiplication maps
$$\cM_U:\bigoplus_{\xi\in U}H^0(F\otimes P_\xi^\vee)\otimes H^0(\H\otimes
P_\xi)
\buildrel{\oplus m_\xi}\over\longrightarrow H^0(F\otimes \H)$$
is surjective for any non-empty Zariski-open set $U\subset \hat X$
containing  $J(F)$.
\end{lemma}

\begin{proof}
The argument follows the proof of \cite{us} Theorem 2.5 (although in fact 
the hypotheses allow us to avoid the use of derived categories).
The statement is equivalent to proving the injectivity of the dual
map (note that $F$ is locally free):
$$\cM_U^\vee:{\rm Ext}^g(\H,F^\vee)\buildrel{\prod m_\xi^\vee}\over\longrightarrow
\prod_{\xi\in U}{\rm Hom}(H^0(\H\otimes
P_\xi),H^g(F^\vee\otimes P_\xi)), $$
where the maps $m_\xi^\vee$ are the co-multiplication maps taking an
extension class $e\in {\rm Ext}^g(\H,F^\vee)$ to its connecting map
$H^0(\H\otimes P_\xi)\rightarrow H^g(F^\vee\otimes
P_\xi)$. 
The index hypotheses allows us to write $\cM_U^\vee$ as the composition of the map
on global sections
$\phi:{\rm Ext}^g(H,F^\vee)\rightarrow {\rm Hom}(\widehat{\H},
\widehat{F^{\vee}})$, followed by the evaluation map 
$$ev_U: H^0({\mathcal H}om(\widehat{\H},
\widehat{F^{\vee}}))\buildrel{\prod ev_\xi}\over\longrightarrow
\prod_{\xi\in U}{\mathcal H}om(\widehat{\H},\widehat{F^{\vee}})(\xi),$$ 
where for a sheaf $\E$, we denote $\E(\xi):=\E\otimes k(\xi)$. 
In addition the hypotheses imply, by Lemma \ref{exts} above,
that the map $\phi$ is an isomorphism. 

On the other hand, if $U$
is a Zariski-open set containing $J(F)$, the map $ev_U$ is injective:
note that by Nakayama's Lemma, given a non-zero global section $s$ of a sheaf
$\E$ on $X$, we have that $s(x)\in
\E(x)$ vanishes identically on a Zariski-open set
$U$ only if either $U$ does not meet a component of the support of the
torsion part of $\E$, or if the torsion part of $\E$ is a sheaf on a non-reduced
subscheme of $X$. Taking $\E={\mathcal H}om(\widehat H,\widehat F^\vee)\cong
(\widehat H)^\vee\otimes \widehat F^\vee$ in our case, 
we see that the torsion part $\tau({\mathcal H}om(\widehat
H,\widehat F^\vee))$ is isomorphic to $(\widehat H)^\vee\otimes\tau(\widehat F^\vee)$. 
The hypothesis that
$U$ contains the Fourier jump locus of $F$ excludes the first possibility
(since the support of the torsion part of $\widehat F^\vee$ is
certainly contained in $J(F)$). The second possibility is excluded by hypothesis. 
\end{proof} 

\begin{remark} More generally, Lemma \ref{mult} and, consequently, Theorem
\ref{WIT} continue to hold (with the same proof) under the hypotheses that
$F$ is a locally free sheaf on an $m$-dimensional Cohen-Macaulay subvariety $Y$
of $X$, ${\mathcal H}om(F,\omega_Y)={\mathcal E}xt^{g-m}(F,\OO_X)$ satisfies W.I.T.
with index $m$, and the torsion part of the Fourier transform ${{\mathcal
H}om(F,\omega_Y)}^{\widehat{}}$ is a torsion-free sheaf on a reduced subvariety of $X$.
\end{remark}

\section{\bf Relative Pontrjagin products}

\subsection*{\bf Pontrjagin products,
multiplication maps and relative Pontrjagin products.}

One of the key points emphasized in
\cite{pareschi} is the relation between multiplication maps of
sheaves on abelian varieties (which are in turn involved in the 
study of linear series) and (skew) Pontrjagin products (cf. Proposition 5.2 below).
Here we develop an analogue relative setting for skew Pontrjagin products,
required for our applications.

\begin{terminology/notation}{\bf ( Skew Pontrjagin product, P.I.T. and $~~~$
P.W.I.T.)}\label{PIT} 
Let us recall first that, given two sheaves $\E$ and $\G$ on $X$, 
their \emph{skew Pontrjagin product} (see \cite{pareschi} \S1) is defined as 
$$\E\hat * \G:={p_1}_*((p_1+p_2)^*(\E)\otimes p_2^*(\G)).$$
We will see, in the spirit of \cite{mukai} \S3, the skew Pontrjagin product as a bifunctor from
${\rm Mod}(X)\times {\rm Mod}(X)$ to
${\rm Mod}(X)$, and we denote by $\buildrel {\bf R}\over{\hat *}$ 
its derived functor. Moreover we adopt the following terminology: 
the pair $(\E,\G)$ \emph{satisfies the Pontrjagin Index Theorem (P.I.T.) 
with index $k=p(\E,\G)$} if 
$h^i((T_x^*\E)\otimes \G)=0$ for any $i\ne k$ and for any $x\in X$.
If $R^i{p_1}_*((p_1+p_2)^*(\E)\otimes {p_2}^*(\G))=0$ for $i\ne
k$ we will say that $(\E,\G)$ \emph{satisfies the Weak Pontrjagin Index Theorem
(P.W.I.T.) with index $k=p(\E,\G)$} and in this case (by abuse of
notation) we will denote again 
$$\E\hat * \G=R^k{p_1}_*((p_1+p_2)^*(\E)\otimes {p_2}^*(\G)).$$
\end{terminology/notation} 

We will also use the following notation: given two sheaves $\E$ and $\G$ on $X$, 
we denote by $\cM(E,G)$ the locus of $x\in X$ where 
the multiplication map 
$$m_x:H^0(T_x^*\E)\otimes H^0(\G)\rightarrow H^0((T_x^*\E)\otimes \G)$$
is not surjective. The relationship between skew Pontrjagin products and 
multiplication maps is provided by the following:

\begin{proposition}\label{mult-pontr}(\cite{pareschi} Proposition 1.1)
Let $\E$ and $\G$ be sheaves on $X$ such that $(\E,\G)$ satisfies P.I.T. with
$p(E,G)=0$. Then $$\cM(\E,\G)=B(\E\hat * \G).$$
\end{proposition}

\begin{remark}
If $E$ and $G$ are locally free 
and $(E,G)$ satisfies P.I.T. with $p(E,G)=0$, then $(E^\vee,G^\vee)$
satisfies P.I.T. with $p(E^\vee,G^\vee)=g$ and, by relative Serre duality,
$(E\hat * G)^\vee\cong
E^\vee\hat * F^\vee$. In other words the dual of
$E\hat * F$ is also a skew Pontrjagin product. 
\end{remark}

In view of Theorems \ref{F-reg} and \ref{WIT}, given a pair
$(E,G)$ satisfying P.I.T. with $p(E,G)=0$ as in the remark above, in order 
to study the surjectivity of the multiplication map
$m_0:H^0(F)\otimes H^0(G)\rightarrow H^0(F\otimes G)$ it is then natural to
investigate whether there exists an ample line bundle
$A$ on $X$ such that the "mixed" product 
$((E\hat * G)\otimes A^\vee)^\vee \cong (E^\vee\hat * F^\vee)\otimes A$
satisfies W.I.T. with $i((E^\vee\hat * F^\vee)\otimes A)=g$. Following
\cite{pareschi}, an appropriate strategy turns out to be
the following: first one establishes a suitable result of "exchange of
Pontrjagin and tensor product under cohomology". Then, to prove the
required vanishing, one uses the fact that, when the sheaves
involved are line bundles algebraically equivalent to powers of a
given one, say $A$, there is a suitable positive integer $n$ such that,
pulling back via multiplication by $n$, the skew Pontrjagin product is
a trivial bundle tensored by a suitable power of $A$. Here we generalize
this technique to a relative setting: Proposition 5.5 below is the
relative analogue of the exchange of Pontrjagin and tensor product under
cohomology while  Proposition 5.6 provides the formulas for the pullback via
multiplication by an integer (both in the usual and the relative
setting). 

\begin{terminology/notation}{\bf (Pontrjagin product relative with respect to
the second variable, relative P.I.T. and P.W.I.T.)}
We denote by $p_i$ and $p_{ij}$ respectively, the projections and
the intermediate projections of
$X\times X\times \hat X$. Consider the bifunctor
$$?~{\hat *}_{rel} ?={p_{13}}_*((p_1+p_2)^*?\otimes p_2^*?\otimes
p_{23}^*{\mathcal P})$$ from ${\rm Mod}(X)\times {\rm Mod}(X)$ to
${\rm Mod}(X\times \hat X)$, and let 
${\buildrel {\bf R}\over{\hat *}}_{rel}$ be its derived functor. 
As usual, we have corresponding notions of Index Theorem and Weak
Index Theorem: e.g. relative P.I.T.
with $p_{rel}(\E,\G)=k$ means that $H^i((T_x^*\E)\otimes \G\otimes
P_\xi)=0$ for any $x\in X$, $\xi \in \hat X$ and $i\ne k$. As
above, if relative P.W.I.T. holds, we write ${\hat *}_{rel}$
rather than ${\buildrel {\bf R}\over{\hat *}}_{rel}$. We denote by $\Gamma$ 
the global sections functor.
\end{terminology/notation}

\begin{proposition}\label{exchange}{\bf (Exchange of Pontrjagin
and tensor product under (relative) cohomology.)}
(a) Assume that $G$ and $H$ are locally free sheaves
on $X$ and $?$ is either an object or a morphism in ${\bf D}(X)$. Then
\newline
\noindent
(i) ${\bf R}\Gamma ((?\buildrel {\bf R}\over{\hat *}G)\otimes
H)\cong {\bf R}\Gamma ((?{\buildrel {\bf R}\over{\hat *}}_{rel} H)
\otimes p_X^*G)$.
\newline
\noindent
(ii) ${\bf R}\hat S((?\buildrel {\bf R}\over{\hat *}G)\otimes
H)\cong {\bf R}{p_{\hat X}}_*((?{\buildrel {\bf R}\over{\hat *}}_{rel} H)
\otimes p_X^*G)$.

\medskip
\noindent
(b) Let in addition $\E$ be a sheaf on $X$ such that
$(\E,G)$ satisfies P.I.T. with $p(\E,G)=k$. Then 
\newline
\noindent
(i) If $(\E,H)$ satisfies P.I.T. with
$p(\E,H)=k$, then  for any $i$,
$$H^i((\E\hat * G)\otimes H)\cong H^i((\E{\hat
*} H)\otimes G).$$
\newline
\noindent
(ii) If $(\E,H)$ satisfies relative P.I.T. with
$p_{rel}(\E,H)=k$, then for any $i$,
$$R^i\hat{\mathcal S}((\E\hat * G)\otimes H)\cong R^i{p_{\hat
X}}_*((\E{\hat *}_{rel} H)\otimes p_X^*G).$$
\end{proposition}

\begin{proof} Part (b)(i) is precisely Lemma 3.2 of \cite{pareschi}. 
The rest of the proof follows the same argument in the derived setting.
Note that we are suppressing part of the symbols showing that we are working with 
the derived functors, in order to simplify the notation.
For the reader's convenience we prove (a)(ii) and (b)(ii). The left hand
side of (a)(ii) can be written as
\begin{eqnarray*}
{\bf R}\hat{\mathcal S}((?\buildrel{\bf R}\over{\hat *} G)\otimes H)& = 
& {\bf R}{p_{\hat X}}_*({\mathcal P}\otimes
p_X^*(H\otimes {\bf R}{p_1}_*((p_1+p_2)^*?\otimes p_2^*G))) \\
& \cong &{\bf R} {p_{\hat X}}_*({\mathcal P}\otimes p_X^*({\bf
R}{p_1}_*(p_1^*H\otimes (p_1+p_2)^*?\otimes p_2^*G))) \\
&\cong & {\bf R}{p_{\hat X}}_*({\bf R}{p_{13}}_*(p_{13}^*\mathcal
P\otimes p_1^*H\otimes (p_1+p_2)^*?\otimes p_2^*G)) \\
&\cong & {\bf R}{p_3}_*
(p_{13}^*{\mathcal
P}\otimes p_1^*H\otimes (p_1+p_2)^*?\otimes p_2^*G) 
\end{eqnarray*}
On the other hand, working out the right hand side we have
\begin{eqnarray*} 
{\bf R}{p_{\hat X}}_*((?{\buildrel {\bf R}\over{\hat *}}_{rel} H)
\otimes p_X^*G) & = &
{\bf R}{p_{\hat X}}_*({\bf R}{p_{13}}_*((p_1+p_2)^*?\otimes
p_2^*H\otimes p_{23}^*{\mathcal P}_{23}) \otimes p_X^*G) \\
&\cong & {\bf R}{p_{\hat X}}_*({\bf
R}{p_{13}}_*(p_1^*G\otimes (p_1+p_2)^*?\otimes p_2^*H\otimes
p_{23}^*{\mathcal P}) \\
& \cong & {\bf R}{p_3}_*(p_1^*G\otimes (p_1+p_2)^*?
\otimes p_2^*H\otimes
p_{23}^*{\mathcal P}) 
\end{eqnarray*}
The result follows via the automorphism $(x,y,\xi)\mapsto
(y,x,\xi)$ of $X\times X\times \hat X$. As for
(b)(ii), under the hypothesis at hand
$(\E\buildrel {\bf R}\over{\hat *}G)\otimes H$ reduces to $(\E\hat *
G)\otimes H[k]$ and $(\E{\buildrel {\bf R}\over{\hat *}}_{rel} H)
\otimes p_X^*G$ reduces to $(\E{\hat
*}_{rel} H)\otimes p_X^*G[k]$. Therefore the assertion follows from (a)(ii).
\end{proof}

\subsection*{\bf Pulling back via multiplication by an integer.}
When line bundles are 
involved, Pontrjagin products usually look simpler when pulled back via
multiplication by an appropriate integer. The purpose
of this subsection is to generalize the results of 
\cite{pareschi} \S3(b) to Pontrjagin products relative with respect
to the second variable. Given a positive integer $n$, the map $x\mapsto
nx$ will be denoted $n_X:X\rightarrow X$.

\begin{proposition}\label{calculations}
(a)  Let $n$ be a
positive integer and let $L$ be a line bundle on~$X$. 
\newline
\noindent
(i) $n_X^*(L \buildrel{{\bf R}}\over {\hat
*} ?)\cong  (L^n \buildrel{{\bf R}}\over {\hat *} (?\otimes
L^{-n+1}))\otimes  n_X^*(L)\otimes p_X^*L^{-n}$
\newline
\noindent
(ii) $(n_X,1_{\hat X})^*(L{\buildrel{{\bf R}}\over {\hat *}}_{rel}?)\cong 
(L^n {\buildrel{{\bf R}}\over {\hat *}} (?\otimes L^{-n+1}))\otimes 
(n_X,1_X)^*(L)\otimes p_X^*L^{-n}$
\newline 
\noindent
(b) Skew Pontrjagin products with the structure sheaf can be expressed as follows:
\newline
\noindent
(i) $?{\buildrel{{\bf R}}\over {\hat *}} {\mathcal O}_X\cong 
{\bf R}\Gamma(?)\otimes \OO_X$
\newline
\noindent
(ii) $?{\buildrel{{\bf R}}\over {\hat *}}_{rel} {\mathcal O}_X\cong
p_{\hat X}^*({\bf R}{\mathcal S}(?))\otimes {\mathcal P}^\vee $
\newline
\noindent
(c) Let $A$ be an ample line bundle on $X$ and assume that
$a$ and $a+b$ are positive integers. Then
\newline
\noindent
(i) $(a+b)_X^*(A^a
 {\hat *} (A^b\otimes P_\xi))\cong H^0(A^{a+b}\otimes P_\xi)\otimes (a+b)_X^*A^a\otimes
a_X^*(A^{-a-b})
\otimes P_\xi^{-a}$
\item[] $(a+b)^*_X(A^{-a}\hat *(A^{-b}\otimes P_\xi))\cong
H^g(A^{-a-b}\otimes P_\xi)\otimes (a+b)_X^*A^{-a}\otimes a_X^*(A^{a+b})\otimes
P_\xi^{a} $
\newline
\noindent
(ii) $((a+b)_X \times 1_{\hat X})^*(A^a
 {\hat *}_{rel} A^b)\cong p_{\hat X}^*((A^{a+b})^{\widehat{}})
\otimes p_X^*((a+b)_X^*A^a\otimes a_X^*A^{-a-b})\otimes 
{\mathcal P}^{-a}$
\item[]
$((a+b)_X \times 1_{\hat X})^*(A^{-a}
{\hat *}_{rel} A^{-b}) \cong  p_{\hat
X}^*((A^{-a-b})^{\widehat{}})
\otimes p_X^*((a+b)_X^*A^{-a}\otimes a_X^*A^{a+b})\otimes
{\mathcal P}^{a}$
\end{proposition}
\begin{proof} Note that (a)(i) is Proposition 3.4 of \cite{pareschi}
and the proof of (a)(ii) identical.
Furthermore, (b)(i) is Remark 3.5(b) of \cite{pareschi} and (b)(ii) is proved in the
same way. Therefore all these proofs are omitted. The first isomorphism of
(c)(i) is Proposition 3.6 of \cite{pareschi}, but note that in that
paper there is a misprint: the last factor of the right hand side of
\cite{pareschi} Prop. 3.6 should read
$a_X^*(A^{-a-b}\otimes \alpha^\vee)$ instead of 
$a_X^*(A^{-a-b})\otimes \alpha^\vee$. The present
formulation follows since $a_X^*P_\xi^\vee=P_\xi^{-a}$. The second
isomorphism of (c)(i) follows by duality. Finally (c)(ii) is proved exactly in
the same way, using Lemma \ref{standard}(a) below, and therefore its proof is omitted too.
\end{proof}

\begin{lemma}\label{standard}
(a)$(n_X, 1_{\hat X})^*{\mathcal P}=
(1_X, n_{\hat X})^*{\mathcal P}={\mathcal P}^{\otimes n}$. 
\newline
\noindent
(b) $R^i{p_{\hat X}}_*({\mathcal P}^{\otimes n})= 0$ for $i< g$ 
and $={\mathcal O}_{{\hat X}_n}$ for $i=g$.
\end{lemma}
\begin{proof} (a) By double duality it is enough to prove
the first equality. We prove it by induction on $n$. For $n=2$, we 
apply the Theorem of the Cube (\cite{mumford} Cor.2 p.58) with 
$f(x,\xi)=g(x,\xi)=(x,0)$, 
$h(x,\xi)=(0,\xi)$ (all maps $X\times \hat X\rightarrow X\times \hat X$)
and $L={\mathcal P}$. We get, using that ${\mathcal P}_{|\{0\}\times \hat X}
={\mathcal O}_{\hat X}$ and ${\mathcal P}_{|X\times \{0\}}={\mathcal O}_X$,
that $(2_X, 1_{\hat X})^*{\mathcal P}={\mathcal P}^{\otimes 2}$. The general formula
follows by induction, applying the same method with 
$f(x,y)=((n-1)_X,0)$, $g(x,y)=(x,0)$, $h(x,y)=(0,y)$.
\newline
\noindent
(b) By flat base change,   
$R^i{p_{\hat X}}_*((1_X\times n_{\hat X})^*({\mathcal P}))
=n_{\hat X}^*R^ip_{\hat X}^*{\mathcal P}= 0$ if $i<g$ and 
$n_{\hat X}^*{\mathcal O}_0={\mathcal O}_{{\hat X}_n}$ if $i=g$.
\end{proof}

\subsection*{\bf A first application: theta-group-free proof of a theorem of Ohbuchi.} 
We end this section by presenting a
theta-group-free proof of (part of) a classical theorem of Ohbuchi
(\cite{ohbuchi2}, see also
\cite{kempf3} Theorem 10.4, \cite{lange} Theorem 7.2.3 and \cite{khaled}
for a proof working in ${\rm char}(k)\ne 2$) on the normal generation of a line
bundle of the form $A^2$, where
$A$ is an ample line bundle on an abelian variety. 
This is intended to be a toy version and an introduction to the new results on 
equations and syzygies in the next section, based on the techniques described above.

Within this framework, it is natural to state Ohbuchi's
Theorem in a slightly more general way.
(This can be however easily deduced from the usual
statement: compare  \cite{lange} Theorem
7.2.3 and Exercise 7.2.) Given an (ample) line bundle
$A$ on $X$, let us denote $s(A):=A\otimes (-1)_X^*A^\vee$.  The map
$s:{\rm Pic}^{c_1(A)}(X)\rightarrow \Pic0$ is surjective and flat
($s(A)\in \Pic0$ classifies the "non-symmetry" of the line bundle $A$). 
Let also $t(A)$ denote a square root of $s(A)$.

\begin{theorem}\label{ohbuchi} ($char(k)\ne 2$) Let $A$ be an ample line
bundle on $X$. Then
$$\mathcal{M}(A^2,A^2)=\bigcup_{\xi\in{\hat X}_2}B(A\otimes P_{t(A)}\otimes P_\xi).$$
Hence $A^2$ is normally generated if and only if $0\not\in \bigcup_{\xi\in
{\hat X}_2}B(A\otimes P_{t(A)}\otimes P_\xi)$.
\end{theorem}

\begin{proof} We will prove only the "positive part" of the result, i.e. the inclusion
$\mathcal{M}(A^2,A^2)\subset\cup_{\xi\in{\hat X}_2}B(A\otimes P_{t(A)}\otimes P_\xi)$,
since this is the part to be generalized in the next section.
The opposite inclusion can be proved similarly.
We have that
\begin{equation}\label{iso}
R^i\hat{\mathcal S}( (A^{-2}\hat *
A^{-2})\otimes A)\cong \left\{\begin{array}{ll}
0 & \textrm{if $i<g$} \\
\widehat{A^{-1}}\otimes{\mathcal O}_{{\hat X}_2+t(A)} &\textrm{if
$i=g$,}
\\ 
\end{array}\right.
\end{equation} 
where $\hat{X}_2+t(A)$ denotes the set $\{~\eta~|~\eta - t(A)\in \hat{X}_2\}$.
Postponing the proof of (\ref{iso}) for a moment, let us show how
it implies the statement. In fact
(\ref{iso}) yields that the hypothesis of
Theorem \ref{WIT} are fulfilled by $(A^2\hat * A^2)\otimes A^\vee$.
Moreover (\ref{iso}) gives also that 
$J((A^2\hat * A^2)\otimes
A^\vee)=\hat{X}_2+t(A)$ (this follows
immediately from base change and Serre duality). Therefore, by
Theorem \ref{WIT}, $B(A^2\hat *A^2)\subset \cup_{\xi\in{\hat X}_2}B(A\otimes
P_{t(A)}\otimes P_\xi)$. On the other hand, by Proposition \ref{mult-pontr},
$\mathcal{M}(A^2,A^2)=B(A^2\hat * A^2)$, and the statement is proved.

\noindent \emph{Proof of (\ref{iso}):} 
\begin{eqnarray}
R^i\hat{\mathcal S}( (A^{-2}\hat *
A^{-2})\otimes A)\cong &&R^i{p_{\hat X}}_*((A^{-2}{\hat
*}_{rel}A)\otimes p_X^*A^{-2}\\
\cong &&R^i{p_{\hat X}}_*(p_{\hat X}^*\widehat{A^{-1}}\otimes
p_X^*(A^{-4}\otimes 2_X^*A)\otimes {\mathcal P}^{2})\\
\cong && R^i{p_{\hat X}}_*(p_{\hat
X}^*\widehat{A^{-1}}\otimes p_X^*(A^{-1}\otimes (-1)_X^*A)\otimes
{\mathcal P}^{2})\\ 
= && R^i{p_{\hat X}}_*(p_{\hat
X}^*\widehat{A^{-1}}\otimes p_X^*P_{s(A)}^\vee\otimes {\mathcal P}^{2}) \\
\cong&&\left\{\begin{array}{ll}
0 & \textrm{if $i<g$} \\
\widehat{A^{-1}}\otimes{\mathcal O}_{{\hat X}_2+t(A)} &\textrm{if $i=g$}.
\\ 
\end{array}\right.
\end{eqnarray}
In the sequence of congruences above, (2) follows by Proposition \ref{exchange}(b)(ii), 
(3) follows by Proposition \ref{calculations}(c)(ii) (second part) with $a=2$ and $b=-1$, 
(4) from the fact that $2_X^*A\cong A^3\otimes (-1)_X^*A$ and (6) from (a slight 
variant of) Lemma 
\ref{standard}(b) and the projection formula.
\end{proof}

\section{\bf Equations defining abelian varieties and
their syzygies}

Putting together the machinery of the previous paragraphs, in 
this section we adress the question of bounding the degrees 
of the generators (and their syzygies) of the homogeneous ideal $I_{X,L}$ of an 
abelian variety $X$ embedded by a complete linear series $|L|$, where
$L$ is a suitable power of an ample line bundle $A$.
Our main result is:

\begin{theorem}\label{main} (${\rm char}(k)\ne 2,3$.)
Let $A$ be an ample line bundle on $X$, with no base divisor. Then: 
\newline
\noindent
(a) If $k\ge 3$ then $I_{X,A^{k}}$ is generated by quadrics.
\newline
\noindent
(b) $I_{X,A^2}$ is generated by quadrics and cubics.
\end{theorem}

Theorem \ref{main} turns out to be a special case of a more general 
result, extending in the now well-known language 
of Green \cite{green} bounds on the degrees of generators of the ideal $I_{X,L}$
to a hierarchy of conditions about higher syzygies. Specifically,
given a variety $X$ embedded in projective space by a complete linear
series $|L|$, the line bundle $L$ is said to \emph{satisfy property
$N_p$ } if the first $p$ steps of the minimal graded free resolution of the
algebra
$R_L=\oplus H^0(L^n)$ over the polynomial ring $S_L=\oplus {\rm Sym}^nH^0(L)$
are linear, i.e. of the form
$$  S_L(-p-1)^{\oplus i_{p}}\rightarrow 
S_L(-p)^{\oplus i_{p-1}}\rightarrow\cdots\rightarrow 
S_L(-2)^{\oplus i_1}\rightarrow S_L\rightarrow R_L\rightarrow 0.$$
Thus $N_0$ means that the embedded variety is projectively normal
(\emph{normal generation} in Mumford's terminology), $N_1$ means that
the homogeneous ideal is generated by quadrics (\emph{normal
presentation}), $N_2$ means that the relations among these quadrics are
generated by \emph{linear} ones and so on. 

More generally even (cf. \cite{pareschi}), one can define properties measuring 
how far the first $p$ steps of the resolution are from being linear. To do this, 
fix $p\ge 0$, and consider the
first $p$ steps of the minimal free resolution of $R_L$ as an $S_L$-module:
$$E_p\rightarrow E_{p-1}\rightarrow \cdots E_1\rightarrow
E_0\rightarrow R_L\rightarrow 0,$$ 
where $E_0=S_L\oplus\bigoplus_jS_L(-a_{0j})$ with $a_{0j}\ge 2$ (since
the linear series is complete),
$E_1=\bigoplus_jS_L(-a_{1j})$ with $a_{1j}\ge 2$ (since 
the embedding is non-degenerate) and so on, up to
$E_p=\bigoplus_jS_L(-a_{pj})$ with $a_{pj}\ge p+1$. Then $L$ is said to 
\emph{satisfy property $N_p^r$} if $a_{pj}\le p+1+r$. In particular 
$N_1^r$ means that $a_{1j}\le 2+r$, i.e. the ideal
$I_{X,L}$ is generated by forms of degree $\le 2+r$, while property 
$N_p^0$ is the same as $N_p$.

\noindent
With this terminology, the extension of Theorem \ref{main} to arbitrary 
syzygies is the following:

\begin{theorem}\label{syzygies} \emph{(${\rm char}(k)$ does not divide $(p+1)$ and
$(p+2)$.)} 
Assume that $A$ has no base divisor. Then:
\newline
\noindent
(a) If $k\ge p+2$ then $A^{k}$ satisfies property $N_p$.
\newline
\noindent
(b) More generally, if $(r+1)(k-1)\ge p+1$ then $A^k$ satisfies
property $N_p^r$.
\end{theorem}

\noindent
A word about the proofs. Although Theorem \ref{main} is subsumed in Theorem 
\ref{syzygies}, we prefer to start by proving it separately. The reason is that 
a substantially higher degree of technicality in the proof of Theorem \ref{syzygies}
might potentially make the main idea less transparent -- with this separation, some 
of the similar arguments will not be repeated in the second proof.

\subsection*{\bf Background material.} We briefly recall some well-known 
facts about the relationship between condition $N_p$, or more generally $N_p^r$,
and the surjectivity of suitable
multiplication maps of vector bundles. For the facts surveyed here see e.g. 
\cite{lazarsfeld} and \cite{pareschi}.
The main point is that condition $N_p^r$ is equivalent
to the exactness in the middle of the piece of the Koszul complex (cf. \cite{green}):
\begin{equation}\label{koszul}
\Lambda^{p+1}H^0(L)\otimes H^0(L^h)\rightarrow \Lambda^p H^0(L)\otimes
H^0(L^{h+1})\rightarrow \Lambda^{p-1}H^0(L)\otimes H^0(L^{h+2})
\end{equation}
for all $h\geq r+1$.
One can in turn express this as a vanishing condition for the cohomology 
of a suitable vector bundle. Specifically, for a globally generated line 
bundle $L$, let $M_L$ be the kernel of the evaluation map:
\begin{equation}\label{ML}
0\rightarrow M_L\rightarrow H^0(L)\otimes {\mathcal O}_X\rightarrow
L\rightarrow 0
\end{equation} 
It follows easily that the exactness of
(\ref{koszul}) is equivalent to the surjectivity of the map
$\Lambda^{p+1}H^0(L)\otimes H^0(L^h)\rightarrow H^0(\Lambda^p M_L\otimes
L^{h+1})$ arising from the exact sequence (obtained by taking exterior powers in (\ref{ML})):
$$0\rightarrow \Lambda^{p+1}M_L\otimes L^h\rightarrow
\Lambda^{p+1}H^0(L)\otimes L^h\rightarrow \Lambda^pM_L\otimes
L^{h+1}\rightarrow 0.$$
Therefore $N_p^r$ holds as soon as 
\begin{equation}\label{wedge}
H^1(\Lambda^{p+1}M_L\otimes L^{h})=0, ~\forall~h\geq r+1.
\end{equation} 
(On abelian varieties the converse is also
true since $H^1(L^h)=0$ for $h\ge 1$.) This leads to:

\begin{proposition}\label{tensor}
(a) \emph{$({\rm char}(k)$ does not divide $(p+1)$)} If $H^1(M_L^{\otimes
(p+1)}\otimes L^h)=0$ for all $h\ge r+1$, then $L$ satisfies condition $N_p^r$.
\newline
\noindent
(b) Assume that $H^1(M_L^{\otimes p}\otimes L^h)=0$. Then
$H^1(M_L^{\otimes (p+1)}\otimes L^h)=0$ if and only if the multiplication map
$$H^0(L)\otimes H^0(M_L^{\otimes p}\otimes L^h)\rightarrow
H^0(M_L^{\otimes p}\otimes L^{h+1})$$ 
is surjective.
\end{proposition}
\begin{proof} Part (a) follows from (\ref{wedge}) since, under the assumption
on the characteristic,
$\Lambda^{p+1}M_L$ is a direct summand of $M_L^{\otimes (p+1)}$. 
Part (b) follows
from the exact sequence
\begin{equation}\label{sequence}
0\rightarrow M_L^{\otimes (p+1)}\otimes L^h\rightarrow H^0(L)\otimes
M_L^{\otimes p}\otimes L^h\rightarrow M_L^{\otimes p}\otimes L^{h+1}\rightarrow 0.
\end{equation}
\end{proof}

\subsection*{\bf Proof of Theorem \ref{main}.} (a) The result is known
for $k\ge 4$, so it is enough to prove it for $k=3$. By Proposition \ref{tensor}(a), 
it suffices to show that
$$H^1(M_{A^3}^{\otimes 2}\otimes A^{3h})=0, ~\forall~h\geq 1.$$
Moreover, we know that $H^1(M_{A^3}\otimes A^{3h})=0$ for
$h\ge 1$ -- by (\ref{sequence}) for $p=0$, this is equivalent to the \emph{normal
generation} of $A^3$, i.e. Koizumi's Theorem.
Therefore, by Proposition \ref{tensor}(b), it is
enough to prove that the multiplication map
$$H^0(A^3)\otimes H^0(M_{A^3}\otimes A^{3h})\rightarrow 
H^0(M_{A^3}\otimes A^{3(h+1)})$$ 
is surjective for $h\ge 1$. Again, this is well known for $h\ge 2$ -- it is
equivalent to the fact that the homogeneous ideal of $X$ embedded by
$|A^3|$ is generated by forms of degree $2$ and $3$ (cf. \cite{kempf1},\cite{lange} 
or, in this interpretation, \cite{pareschi}). Therefore the only case to be examined is
$h=1$.  

We prove more generally that the locus
$\mathcal{M}(A^3,M_{A^3}\otimes A^3)$ is empty (cf. Proposition \ref{mult-pontr}). 
By Proposition \ref{mult-pontr} we have 
$\mathcal{M}(A^3,M_{A^3\otimes A^3})=  B(A^3\hat * (M_{A^3}\otimes A^3))$,
since from the defining sequence \ref{ML} it is not hard to see that the pair
$(A^3,M_{A^3}\otimes {A^3})$ satisfies P.I.T. with index $0$. 
To this end we make the following:

\noindent
\emph{Claim.
If $A$ has no
base divisor, then $(A^3\hat * (M_{A^3}\otimes A^3))\otimes A^\vee$ is 
$M$-regular.}

\noindent  
By Theorem \ref{F-reg} this yields that  $A^3\hat *
(M_{A^3}\otimes A^3)$ is globally generated, and hence the theorem. 

\noindent
\emph{Proof of Claim.} 
Recall from \cite{us} \S3 that it is enough to prove that the cohomological 
support loci
$$V^i:=\{~\xi\in \hat X\> | \> h^i((A^3\hat *( M_{A^3}\otimes A^3))\otimes
A^\vee\otimes P_\xi)>0\}$$
have {\rm codim}ension $> i$ for all $i>0$. By Proposition \ref{exchange}(b)(i) we
have that
$$V^i=\{\xi\in \hat X \> |\> h^i((A^3\hat *(A^\vee \otimes P_\xi))\otimes
M_{A^3}\otimes A^3)>0\}.$$
Let us consider the exact sequence obtained from \ref{ML}
\begin{eqnarray*}
&0\rightarrow (A^3\hat *(A^\vee \otimes P_\xi))\otimes
M_{A^3}\otimes A^3\rightarrow H^0(A^3)\otimes 
(A^3\hat *(A^\vee \otimes P_\xi)) \otimes A^3\rightarrow&\\
&\rightarrow (A^3\hat *(A^\vee \otimes P_\xi)) \otimes A^6\rightarrow 0&
\end{eqnarray*}

\vskip0.2truecm\noindent\emph{Subclaim. 
$h^i(((A^3\hat *(A^\vee \otimes P_\xi)) \otimes A^n)=
0$ for any $n\ge 3 $, $\xi\in \hat X$ and $i\ge 1$.}

\begin{proof} This is again a standard application
of Proposition \ref{calculations}(c)(i): taking $a=3$ and $b=-1$ we have
that $2_X^*(A^3\hat * (A^\vee\otimes P_\xi))\cong H^0(A^2\otimes P_\xi)\otimes
2_X^*A^3\otimes 3_X^*A^{-2}\otimes P_\xi^{-3}$. This implies that 
$$2_X^*(A^3\hat * (A^\vee\otimes P_\xi)\otimes A^n)\cong
H^0(A^2\otimes P_\xi)\otimes
2_X^*A^3\otimes 3_X^*A^{-2}\otimes P_\xi^{-3}\otimes 2_X^*A^n,$$ 
which is isomorphic to a sum of copies of line bundles algebraically equivalent
to $A^{(4n-6)}$, thus certainly ample for $n\ge 3$. As we are in
characteristic $\ne 2$, $H^i((A^3\hat *(A^\vee \otimes P_\xi)) \otimes A^n)$
is a direct summand of $H^i(2_X^*((A^3\hat *(A^\vee \otimes P_\xi))\otimes
A^n))$, which proves the subclaim.
\end{proof}
 
\noindent
Passing to cohomology in the exact sequence above, by the Subclaim  we have that
\begin{enumerate}
\item[(i)]$V^i$ is empty for $i\ge 2$. 
\item[(ii)] $V^1$ coincides with the locus where the
multiplication map
$$H^0(A^3)\otimes H^0(
(A^3\hat *(A^\vee \otimes P_\xi)) \otimes A^3)\rightarrow
H^0((A^3\hat *(A^\vee \otimes P_\xi)) \otimes A^6)$$
is not surjective.
\end{enumerate}

\noindent 
In view of (i), the Claim would be implied by the inequality
 ${\rm {\rm codim}}(V^1)>1$. 
Again by Proposition \ref{mult-pontr}, we have that 
\begin{equation}\label{S1}
V^1=\bigl\{\xi\in \hat X \>|\> 0_X\in
B\bigl(A^3\hat * ((A^3\hat *(A^\vee
\otimes P_\xi))\otimes A^3)\bigr)\bigr\}.
\end{equation}
We will approach this by the same trick of twisting with $A^\vee$ 
in order to try and apply Theorem \ref{WIT}. By relative duality we have
$$((A^3\hat * (A^3\hat *(A^\vee
\otimes P_\xi))\otimes A^3)\otimes A^\vee)^\vee\cong
(A^{-3}\hat * (A^{-3}\hat *(A\otimes P_\xi^\vee))\otimes A^{-3})\otimes
A.$$ By Proposition \ref{exchange}(b)(ii) 
\begin{equation}\label{strunz}
R^i\hat{\mathcal S}((A^{-3}\hat *
(A^{-3}\hat *(A\otimes P_\xi^\vee))\otimes A^{-3})\otimes A)\cong
{R^ip_{\hat X}}_*((A^{-3}{\hat *}_{rel}A)\otimes p_X^*((A^{-3}\hat *
(A\otimes P_\xi^\vee))\otimes A^{-3})).
\end{equation} 
The key point that W.I.T. with index $g$ is satisfied goes through the following:
\begin{eqnarray}\label{iso2}
&R^i{p_{\hat X}}_*\bigl((2_X,1_{\hat X})^*[(A^{-3}{\hat
*}_{rel}A)\otimes p_X^*((A^{-3}\hat * (A\otimes P_\xi^\vee))\otimes A^{-3})]
\bigr)\cong&\nonumber \\
&\cong \left\{\begin{array}{ll}
  0 &\textrm{if $i<g$} \\
H^g(A^{-2}\otimes P_\xi^\vee)\otimes
\widehat{A^{-2}}\otimes{\mathcal
O}_{{\hat X}_3 -s(A)-\xi}
 &\textrm{if
$i=g$}
\end{array}\right.&
\end{eqnarray}
\noindent
\emph{Proof of (\ref{iso2}).}   By Proposition \ref{calculations}(c)(ii) with 
$a=3$ and $b=-1$:
$$(2_X,1_{\hat X})_*( A^{-3}{\hat *}_{rel}A)\cong p_{\hat
X}^*\widehat{A^{-2}}\otimes p_X^*(2_X^*A^{-3}\otimes 3_X^*A^2)\otimes
{\mathcal P}^{3}.$$ 
In conclusion, using also Proposition \ref{calculations}(c)(i),
\begin{eqnarray*}
&(2_X,1_{\hat X})^*[(A^{-3}{\hat
*}_{rel}A)\otimes p_X^*((A^{-3}\hat * (A\otimes P_\xi^\vee))\otimes A^{-3})]\cong
&\\
 & p_{\hat X}^*\widehat{A^{-2}}\otimes
p_X^*(2_X^*A^{-3}\otimes 3_X^*A^2\otimes H^g(A^{-2}\otimes P_\xi^\vee)\otimes
2_X^*A^{-3}\otimes 3_X^*A^2 \otimes P_\xi^{-3}\otimes 2_X^*A^{-3})\otimes {\mathcal P}^{3}
\cong & \\
&\cong p_X^*H^g(A^{-2}\otimes P_\xi^\vee)\otimes p_{\hat X}^*\widehat{A^{-2}}\otimes
p_X^*(P_{s(A)}^\vee\otimes P_\xi^\vee)^3\otimes {\mathcal P}^{3}
\end{eqnarray*}
(Cf. the notation introduced before Theorem \ref{ohbuchi}.) 
Therefore (\ref{iso2}) follows from the projection formula and Lemma
\ref{standard}.

We are now able to conclude the proof of the Claim. As we are in ${\rm char}(k)\ne
2$,  
$$R^i{p_{\hat
X}}_*((A^{-3}{\hat *}_{rel}A)\otimes {p_X}^*((A^{-3}\hat * (A\otimes
P_\xi^\vee))\otimes A^{-3}))$$ 
is a direct summand of 
$$R^i{p_{\hat X}}_*\bigl((2_X,1_{\hat X})^*[(A^{-3}{\hat
*}_{rel}A)\otimes p_X^*((A^{-3}\hat * (A\otimes P_\xi^\vee))\otimes A^{-3})]
\bigr).$$
According to 
(\ref{strunz}) and (\ref{iso2}), it follows that
$(A^3\hat * (A^3\hat *(A^\vee
\otimes P_\xi)\otimes A^3))\otimes A^\vee$ satisfies both hypotheses
of Theorem \ref{WIT}. Moreover, by (\ref{iso2}) it also follows, 
using relative duality and base change, that there is an inclusion 
$$J((A^3\hat * (A^3\hat *(A^\vee
\otimes P_\xi)\otimes A^3))\otimes A^\vee)\subset {\hat
X}_3-s(A)-\xi.$$ 
Therefore Theorem \ref{WIT} implies that  
\begin{equation}\label{inclusion}
B\bigl(A^3\hat * ((A^3\hat *(A^\vee
\otimes P_\xi))\otimes A^3)\bigr) \subset \bigcup_{\eta\in {\hat
X}_3}B(A\otimes P_{-s(A)-\xi}\otimes P_\eta).
\end{equation}
From (\ref{S1}) and (\ref{inclusion}) it follows that
${\rm codim}(S^1)={\rm codim}(B((A^3\hat *(A^\vee
\otimes P_\xi))\otimes A^3))\ge {\rm codim}(B(A))$ and this
proves the Claim. 

\medskip
\noindent
(b) The proof is completely similar. This time one has to prove
that the multiplication map $$H^0(A^2)\otimes H^0(M_{A^2}\otimes
A^4)\rightarrow H^0(M_{A^2}\otimes A^6)$$ 
is surjective. As above, the
result follows from the following statement, proved in the same way:
\newline
\noindent
\emph{Claim. If $A$ has no base divisor then
$(A^2\hat * (M_{A^2}\otimes A^4))\otimes A^\vee$ is
$M$-regular.}

\subsection*{\bf Proof of Theorem \ref{syzygies}.} The argument is
a combination between the proof of Lazarsfeld's conjecture in 
\cite{pareschi} and the idea of the
proof of Theorem \ref{main}.  
We prove only (a), since the proof of (b) is completely analogous.
First of all the theorem is known for $k\ge p+3$
and so we need to prove it only for $k=p+2$. For
$L=A^{p+2}$, the exactness of the complex (\ref{koszul}) is known to hold
for $h\ge 2$. 
(This means that the syzygies at the $p$-th
step are generated at most in degree 2, i.e. condition $N_p^1$ in
terminology of
\cite{pareschi}, and so it follows from \cite{pareschi} Theorem 4.3.)
Putting everything together, by Proposition \ref{tensor} it follows that it 
suffices to prove that the multiplication map
$$H^0(A^{p+2})\otimes H^0(M_{A^p+2}^{\otimes p}\otimes A^p)\rightarrow
H^0(M_{A^p+2}^{\otimes p}\otimes A^{2p}) $$
is surjective. 
Given $\xi\in\hat X$, we denote
$$F_\xi^{(n,m)}:=A^n\hat * (A^m\otimes P_\xi)$$ 
We will prove the following:

\vskip0.2truecm\noindent\emph{Claim 1.}  \emph{For every integer $k$ such
that $1\le k\le p$  and every $\xi_1,\dots ,\xi_k\in \hat X$,
the locally free sheaf
$$A^{p+2}\hat * (M_{A^{p+2}}^{\otimes
k} \otimes\bigotimes_{i=1}^{p-k} F_{\xi_i}^{(p+2,-1)}
\otimes A^{p+2})$$ is globally generated.}

\vskip0.2truecm\noindent For $k=p$ this, together with Proposition
\ref{mult-pontr}, proves the theorem (again, the fact that P.I.T. with index
$0$ is verified follows from \cite{pareschi} Proposition 4.2). 

\vskip0.2truecm\noindent\emph{Proof of Claim 1.} This goes by induction on $k$.
Let us assume for a moment that we know the initial step $k=1$, and show that 
the statement for $k-1$ implies the statement for $k$, for all $k\geq 2$.

We fix $\xi_1, \dots ,\xi_k\in \hat X$.  By Theorem
\ref{F-reg}, it is enough to prove that the vector bundle
$$(A^{p+2}\hat * (M_{A^{p+2}}^{\otimes
k} \otimes\bigotimes_{i=1}^{p-k} F_{\xi_i}^{(p+2,-1)}
\otimes A^{p+2}))\otimes A^\vee$$
is $M$-regular.
In fact, for $k\geq 2$, it will even satisfy I.T. with index $0$,
i.e. 
\begin{equation}\label{mah}
H^i((A^{p+2}\hat * (M_{A^{p+2}}^{\otimes
k} \otimes\bigotimes_{i=1}^{p-k} F_{\xi_i}^{(p+2,-1)}
\otimes A^{p+2}))\otimes A^\vee\otimes P_\xi)=0
\end{equation}
for any $i>0$ and any $\xi\in \hat X$. By Proposition
\ref{exchange}(b)(i) we have that
\begin{eqnarray*}
&H^i((A^{p+2}\hat * (M_{A^{p+2}}^{\otimes
k} \otimes\bigotimes_{i=1}^{p-k} F_{\xi_i}^{(p+2,-1)}
\otimes A^{p+2}))\otimes A^\vee\otimes P_\xi)\cong &\\
&\cong H^i((A^{p+2}\hat * (A^\vee\otimes P_\xi))\otimes
 M_{A^{p+2}}^{\otimes
k}\otimes\bigotimes_{i=1}^{p-k} F_{\xi_i}^{(p+2,-1)}
\otimes A^{p+2}).&
\end{eqnarray*}
(As usual, the hypotheses of Proposition \ref{exchange}(b)(ii) are fulfilled
because of \cite{pareschi} Proposition 4.2.) Then, as in the proof of the previous 
theorem (and we won't go through all the details), the
sequence
$$0\rightarrow M_{A^{p+2}}^{\otimes k}\rightarrow H^0(A^{p+2})\otimes
M_{A^{p+2}}^{\otimes k-1}\rightarrow
M_{A^{p+2}}^{\otimes  k-1}\otimes A^{p+2}\rightarrow 0$$ 
twisted by
$(A^{p+2}\hat * (A^\vee\otimes P_\xi))\otimes \bigotimes_{i=1}^{p-k}
F_{\xi_i}^{(p+2,-1)}\otimes A^{p+2}$
gives that the cohomology groups of (\ref{mah}) are zero, except for
$H^1$ which vanishes if only if the multiplication map 
\begin{eqnarray*}
&H^0(A^{p+2})\otimes H^0(M_{A^{p+2}}^{\otimes k-1}\otimes (A^{p+2}\hat *
(A^\vee\otimes P_\xi))\otimes
\bigotimes_{i=1}^{p-k} F_{\xi_i}^{(p+2,-1)}\otimes A^{p+2})\rightarrow&\\
&\rightarrow  H^0(M_{A^{p+2}}^{\otimes k-1}\otimes (A^{p+2}\hat *
(A^\vee\otimes P_\xi))\otimes
\bigotimes_{i=1}^{p-k} F_{\xi_i}^{(p+2,-1)}\otimes A^{2(p+2)})
\end{eqnarray*} is
surjective. But this follows from the inductive
hypothesis and Proposition \ref{mult-pontr}.

We are left with proving Claim 1 for $k=1$. To this end we will apply the same 
reasoning as in the previous paragraph, only this time we use Theorem 
\ref{F-reg} with the weaker input that the sheaf in question is just $M$-regular.

\noindent
\emph{Claim 2. The sheaf $A^{p+2}\hat * (M_{A^{p+2}}
\otimes\bigotimes_{i=1}^{p-1} F_{\xi_i}^{(p+2,-1)}
\otimes A^{p+2})$ is $M$-regular.} 

\noindent
\emph{Proof.} As before, the claim follows if we show that the locus of $\xi$
such that the multiplication map
\begin{eqnarray*}
&H^0(A^{p+2})\otimes H^0((A^{p+2}\hat *
(A^\vee\otimes P_\xi))\otimes
\bigotimes_{i=1}^{p-1} F_{\xi_i}^{(p+2,-1)}\otimes A^{p+2})\rightarrow&\\
&\rightarrow  H^0((A^{p+2}\hat *
(A^\vee\otimes P_\xi))\otimes
\bigotimes_{i=1}^{p-1} F_{\xi_i}^{(p+2,-1)}\otimes A^{2(p+2)})
\end{eqnarray*}
is surjective has codimension at least $2$. By Proposition \ref{mult-pontr}, this 
locus is precisely 
$$\{~\xi~|~0_X\in B(A^{p+2}\hat * (F_\xi^{(p+2,-1)}\otimes
\bigotimes_{i=1}^{p-1} F_{\xi_i}^{(p+2,-1)}\otimes A^{p+2}))\}.$$
By relative duality the dual of 
$$(A^{p+2}\hat * (F_\xi^{(p+2,-1)}\otimes
\bigotimes_{i=1}^{p-1} F_{\xi_i}^{(p+2,-1)}
\otimes A^{p+2}))\otimes A^\vee $$
is 
$$(A^{-p-2}\hat *(
F_{-\xi}^{(-p-2,1)}\otimes
\bigotimes_{i=1}^{p-1} F_{-\xi_i}^{(-p-2,1)}
\otimes A^{-p-2}))\otimes A$$
By Proposition \ref{exchange} we have that
\begin{eqnarray}\label{cagata}
&R^i\hat{\mathcal S}
((A^{-p-2}\hat *
(F_{-\xi}^{(-p-2,1)}\otimes
\bigotimes_{i=1}^{p-1} F_{-\xi_i}^{(-p-2,1)}
\otimes A^{-p-2}))\otimes A)\cong &\nonumber \\
&R^i{p_{\hat X}}_*((A^{-p-2}{\hat *}_{rel}
A)\otimes p_X^*(F_{-\xi}^{(-p-2,1)}\otimes
\bigotimes_{i=1}^{p-1} F_{-\xi_i}^{(-p-2,1)}
\otimes A^{-p-2}))&
\end{eqnarray}
The key point, analogous to (\ref{iso2}) of the previous proof, is that 
\begin{eqnarray}\label{iso3}
&R^i{p_{\hat X}}_*(((p+1)_X,1_{\hat X})^*((A^{-p-2}{\hat *}_{rel}
A)\otimes  F_{-\xi}^{(-p-2,1)}\otimes
\bigotimes_{i=1}^{p-1} F_{-\xi_i}^{(-p-2,1)}
\otimes A^{-p-2}))\cong&\nonumber\\
& \cong \left\{\begin{array}{ll}
  0 &\textrm{if $i<g$} \\
V\otimes
\widehat{A^{-p-1}}\otimes{\mathcal
O}_{{\hat X}_{p+2}-\xi-\sum_{i=1}^{p-1}\xi_i-\frac{(p+1)(p+2)}{2}s(A)}
 &\textrm{if
$i=g$}
\end{array}\right.&
\end{eqnarray}
where $V$ is a suitable vector space.
\newline
\noindent
\emph{Proof of (\ref{iso3}).} By Proposition
\ref{calculations}(c)(i) and (ii) we have that
\begin{eqnarray*}
&((p+1)_X,1_{\hat X})^*((A^{-p-2}{\hat *}_{rel}
A))\otimes F_{-\xi}^{(-p-2,1)}\otimes
\bigotimes_{i=1}^{p-1} F_{-\xi_i}^{(-p-2,1)}
\otimes A^{-p-2})\cong &\\
&p_{\hat
X}^*\widehat{A^{-p-1}}\otimes
V\otimes p_X^*((p+1)_X^*A^{-(p+2)(p+1)}\otimes (p+2)_X^*
A^{(p+1)^2}\otimes P_{-\xi}^{p+2}\otimes\bigotimes_{i=1}^{p-1}P_{-\xi_i}^{p+2})\otimes {\mathcal
P}^{-p-2} &\\
\end{eqnarray*}
where $V$ is the vector space
$\bigotimes_{i=1}^{p-1}H^g(A^{-p-1}\otimes P_{{\xi}_i}^\vee)\otimes 
H^g(A^{-p-1}\otimes P_{\xi}^\vee)$. Therefore
(\ref{iso3}) follows from Lemma \ref{standard}, noting that, by a standard
calculation, $$(p+1)_X^*(A^{-(p+2)(p+1)}\otimes A^{-p-2})\otimes (p+2)_X^*
A^{(p+1)^2}
\cong P_{-s(A)}^{(p+1)^2(p+2)/2}.$$

\noindent
The argument goes now as in the previous proof: we have first that
$$R^i{p_{\hat X}}_*((A^{-p-2}{\hat *}_{rel}
A) \otimes F_{-\xi}^{(-p-2,1)}\otimes
\bigotimes_{i=1}^{p-1} F_{-\xi_i}^{(-p-2,1)}
\otimes A^{-p-2})$$ is a direct summand of 
$$R^i{p_{\hat X}}_*((2_X,1_{\hat X})^*((A^{-p-2}{\hat *}_{rel}
A) \otimes F_{-\xi}^{(-p-2,1)}\otimes
\bigotimes_{i=1}^{p-1} F_{-\xi_i}^{(-p-2,1)}
\otimes A^{-p-2})).$$ 

\noindent
Summing up, the sheaf 
$$(A^{p+2}\hat * (F_\xi^{(p+2,-1)}\otimes
\bigotimes_{i=1}^{p-1} F_{\xi_i}^{(p+2,-1)}
\otimes A^{p+2}))\otimes A^\vee $$
satisfies the hypotheses of Theorem \ref{WIT}, and its Fourier jump locus 
is included in ${\hat X}_{p+2}-\xi-\sum_{i=1}^{p-1}\xi_i-\frac{(p+1)(p+2)}{2}s(A)$.
Thus, by Theorem \ref{WIT}, we finally have that
$$B(A^{p+2}\hat * (F_\xi^{(p+2,-1)}\otimes
\bigotimes_{i=1}^{p-1} F_{\xi_i}^{(p+2,-1)}\otimes A^{p+2}))\subset
\bigcup_{\eta\in {\hat
X}_{p+2}}B(A\otimes P_{-\frac{(p+1)(p+2)}{2}s(A)-\sum\xi_i-\xi}\otimes P_\eta),$$
and the Claim follows since the base locus of $A$ is of codimension at least $2$.

\subsection*{\bf A conjecture based on the $M$-regularity index.}

As already mentioned in Section 3, Theorem \ref{syzygies} raises a natural 
question about a potentially general relationship between the equations and 
syzygies of $X$ in the embedding given by a power of a line bundle $A$, and 
the higher order properties of $A$, reflected in the $M$-regularity index $m(A)$
defined in \S3.

\begin{conjecture}
\emph{  Let $p\ge m$ be non-negative integers.
If $A$ is ample and $m(A)\ge m$, then $A^{\otimes k}$ satisfies
$N_p$ for any $k\ge p+3-m$.  } 
\end{conjecture}

This conjecture is a refinement of Lazarsfeld's conjecture, proved in \cite{pareschi}, 
which is the case $m(A)=0$, i.e. no conditions on $A$. Theorem \ref{syzygies} gives 
an affirmative answer to the conjecture for $m(A)=1$, which by Example \ref{base_div} 
happens precisely when $A$ has no base divisor. We remark though that the methods of 
this paper fail to apply for powers $A^k$ with $k\leq p+1$, so a new idea seems to 
be needed for the case of higher regularity indices.

\providecommand{\bysame}{\leavevmode\hbox to3em{\hrulefill}\thinspace}

\end{document}